\documentclass{elsart}
\usepackage{amssymb,amsfonts}
\newcommand{\trace}{\mathop{\rm Tr}\nolimits}
\newcommand{\tr}{\mathop{\rm tr}\nolimits}

\newcommand{\diag}{\mathop{\rm Diag}\nolimits}

\newcommand{\twomat}[4]{\left(\begin{array}{cc}#1&#2\\#3&#4\end{array}\right)}
\newcommand{\cP}{{\mathcal P}}

\newcommand{\C}{{\mathbb{C}}}
\newcommand{\R}{{\mathbb{R}}}
\newcommand{\N}{{\mathbb{N}}}
\newcommand{\M}{{\mathbb{M}}}
\DeclareRobustCommand\openone{\leavevmode\hbox{\small1\normalsize\kern-.33em1}}
\newcommand{\id}{\mathbf{I}}

\newcommand{\ddt}{\frac{\partial}{\partial t}}
\newcommand{\ddtp}{\left.\frac{\partial}{\partial t}\right|_{t\to 0^+}}

\newcommand{\be}{\begin{equation}}
\newcommand{\ee}{\end{equation}}
\newcommand{\bea}{\begin{eqnarray}}
\newcommand{\eea}{\end{eqnarray}}
\newcommand{\beas}{\begin{eqnarray*}}
\newcommand{\eeas}{\end{eqnarray*}}
\newtheorem{definition}{Definition}
\newtheorem{theorem}{Theorem}
\newtheorem{lemma}{Lemma}
\newtheorem{corollary}{Corollary}

\newtheorem{proposition}{Proposition}
\newtheorem{question}{Question}

\begin{document}
\begin{frontmatter}
\title{On norm sub-additivity and super-additivity inequalities for concave and convex functions}
\author{Koenraad M.R.\ Audenaert}
\address{
Department of Mathematics, Royal Holloway, University of London,\\
Egham TW20 0EX, United Kingdom}
\ead{koenraad.audenaert@rhul.ac.uk}
\author{Jaspal Singh Aujla}
\address{Department of  Mathematics, National Institute of Technology,\\
Jalandhar 144011, Punjab, India}
\ead{aujlajs@nitj.ac.in}

%------------------------------------------------------------------ ABSTRACT
%-----------------------
\begin{abstract}
Sub-additive and super-additive inequalities for concave and convex functions have been generalized
to the case of matrices by several authors over a period of time. These lead to some interesting
inequalities for matrices, which in some cases coincide with, and in other cases are at variance with
the corresponding inequalities for real numbers. We survey some of these matrix inequalities and do
further investigations into these.

We introduce the novel notion of dominated majorization between the spectra
of two Hermitian matrices $B$ and $C$, dominated by a third Hermitian matrix $A$.
Based on an explicit formula for the gradient of the sum of the $k$ largest eigenvalues of a
Hermitian matrix, we show
that under certain conditions dominated majorization
reduces to a linear majorization-like relation between the diagonal elements of $B$ and $C$ in a certain basis.
We use this notion as a tool to give new, elementary proofs for the sub-additivity
inequality for non-negative concave functions first proved by Bourin and Uchiyama
and the corresponding super-additivity inequality for non-negative convex functions
first proven by Kosem.

Finally, we present counterexamples to some conjectures that
Ando's inequality for operator convex functions could more generally hold,
e.g.\ for ordinary convex, non-negative functions.

\bigskip
\begin{center}
\textit{Dedicated to the memory of Ky Fan}
\end{center}
\end{abstract}

\end{frontmatter}
%%%%%%%%%%%%%%%%%%%%%%%%%%%%%%%%%%%%%%%%%%%%%%%%%%%%%%%%%%%%%%%%%%%%%%%%%%%%%%%%%%%%%%%%%%%%%%%
%%%%%%%%%%%%%%%%%%%%%%%%%%%%%%%%%%%%%%%%%%%%%%%%%%%%%%%%%%%%%%%%%%%%%%%%%%%%%%%%%%%%%%%%%%%%%%%
\section{Introduction}
Two of the basic properties that a real-valued function $f(x)$ defined over the reals can possess
are sub-additivity and super-additivity. Sub-additivity means that for all $x$, $y$ in the domain of $f$,
$$
f(x+y) \le f(x)+f(y),
$$
while super-additivity means the opposite
$$
f(x)+f(y) \le f(x+y).
$$
Two classical theorems that characterise sub- and super-additivity for functions defined on $\R^+$
(although not completely)
are presented as Theorem 7.2.4 and 7.2.5 in \cite{hillephillips}.
Their Theorem 7.2.4 states that functions $f$ for which $f(t)/t$ is decreasing
in $\R_+$ are subadditive.
Theorem 7.2.5 in \cite{hillephillips} states that any measurable concave function $f$ is subadditive in $\R_+$
iff $f(0+)\ge 0$.

In recent years, ongoing effort has been spent to characterise matrix functions exhibiting
similar sub-additivity or super-additivity properties.
Of course, many variations on this theme are possible, and in this paper we restrict attention
to sub- and super-additivity in norm for non-negative functions.
For a given unitarily invariant norm $|||\cdot|||$, these amount
to the norm inequalities
$|||f(A)+f(B)|||\le |||f(A+B)|||$ (or reversed), with positive semidefinite $A$ and $B$,
but one can equally well consider the inequality
$|||f(A)-f(B)|||\le |||f(|A-B|)|||$ (or reversed).
Historically, these inequalities have been proven first for operator monotone, and/or operator concave
functions $f$, and only later have they been generalised to non-negative functions that are concave and/or convex.
Interestingly, the proofs of these generalisations exploit the corresponding results for
operator monotone/concave functions.

In this paper we first give a historical overview of these developments, in Sections \ref{sec:main}
and \ref{sec:mainmin}. Then we resolve a number of
still open questions regarding the inequality $|||f(|A-B|)|||\le |||f(A)-f(B)|||$,
which is known to be true for operator convex functions. We show by counterexample
that it does not hold in general for non-negative convex functions, nor do a number of successively weakened
versions. By imposing the condition $A\ge ||B||_\infty$, we obtain the closest match of this
inequality that does hold for convex functions (or in reversed sense for concave functions), namely
the eigenvalue inequality $\lambda^\downarrow_k(f(A-B)) \le \lambda^\downarrow_k(f(A)-f(B))$, for all $k$.

In Section \ref{sec:apps}, we present
a new and elementary proof of a sub-additivity norm inequality for
non-negative concave functions and a super-additivity norm inequality for non-negative convex functions,
that do not rely on the corresponding inequality for operator monotone/convex functions,
nor on the theory of operator monotone functions.
The proof exploits the novel notion of dominated majorization between the spectra
of two Hermitian matrices $B$ and $C$, dominated by a third Hermitian matrix $A$.
Based on an explicit formula for the gradient of the sum of the $k$ largest eigenvalues of a
Hermitian matrix, we show
that under certain conditions this dominated majorization
reduces to a linear majorization-like relation between the diagonal elements of $B$ and $C$ in a certain basis.
This is explained in full detail in Section \ref{sec:domi} (with one of the proofs postponed to Section
\ref{sec:proofdelta}).
%%%%%%%%%%%%%%%%%%%%%%%%%%%%%%%%%%%%%%%%%%%%%%%%%%%%%%%%%%%%%%%%%%%%%%%%%%%%%%%%%%%%%%%%%%%%%%%
%%%%%%%%%%%%%%%%%%%%%%%%%%%%%%%%%%%%%%%%%%%%%%%%%%%%%%%%%%%%%%%%%%%%%%%%%%%%%%%%%%%%%%%%%%%%%%%
\section{Preliminaries}

In this section, we introduce the notations and necessary prerequisites;
a more detailed exposition can be found, e.g.\ in \cite{bhatia97}.

Throughout, $\M_{n}$ shall denote the set of $n\times n$ complex matrices
and $\M^H_{n}$ shall denote the set of all Hermitian matrices
in $\M_{n}$. We shall abbreviate the terms positive
semidefinite and positive definite by PSD and PD, respectively. By $A\ge B,$ we mean that $A-B\ge0$.
Let $I$ be an interval in $\R$. We shall denote
by $\M^H_{n}(I),$ the set of all Hermitian matrices in $\M_{n}$ whose spectrum is contained in the interval $I$.

We denote the identity matrix by $\id$, and use the shorthand $a=a\id$ for scalar matrices.

We denote the absolute value by $|\cdot|$, both for scalars and for matrices. For matrices this is defined
as $|A| := (A^* A)^{1/2}$. Similarly, we denote the positive part of a real scalar or Hermitian
matrix by $(\cdot)_+$, and define it by $A_+ := (A+|A|)/2$. We denote the vector of diagonal entries of a
matrix $A$ by $\diag(A)$.
We will use the abbreviations LHS and RHS for left-hand side and right-hand side, respectively.

Let $A\in \M^H_{n}(I)$ have the spectral decomposition
\begin{center}
$A=U^{*}{\rm diag}(\lambda _{1},\lambda _{2},\ldots ,\lambda _{n})U$
\end{center}
where $U$ is a unitary matrix and $\lambda _{1},\lambda _{2},\ldots ,\lambda _{n}$ are the eigenvalues of $A$.
Let $f$ be a real
valued function defined on $I$. Then
$f(A)$ is defined by
\begin{center}
$f(A)=U^{*}{\rm diag }(f(\lambda _{1}),f(\lambda _{2}),\ldots ,f(\lambda _{n}))U$.
\end{center}

Let $n\in \N$ be arbitrary but fixed.
The function $f$ is called \textit{matrix monotone of order $n$ on $I$} if
$$
A\ge B ~~\Longrightarrow ~~~ f(A)\ge f(B)
$$
for all $A,B\in \M^H_n(I)$,
and \textit{matrix convex of order $n$ on $I$} if
\begin{center}
$f(\alpha A+(1-\alpha )B)\leq \alpha f(A)+(1-\alpha )f(B)$
\end{center}
for all $0\leq \alpha \leq 1$ and $A,B\in \M^H_n( I)$.
Likewise, $f$ is called \textit{matrix concave of order $n$ on $I$} if
$-f$ is matrix convex of order $n$ on $I$. If the function $f$ is matrix monotone
of all orders $n$ on $I$ then $f$ is called \textit{operator monotone on $I$}.
The operator convexity and operator concavity are defined similarly.

A norm $||| \cdot  |||$ on $\M_{n}$ is called
unitarily invariant (UI) or symmetric if
\begin{center}
$||| UAV ||| ~=~ |||A|||$
\end{center}
for all $A\in \M_{n}$ and for all unitary $U,V\in \M_{n}$.
The most basic unitarily invariant norms are the
Ky Fan norms $|| \cdot || _{(k)}, (k=1,2,\cdots ,n),$
defined as
$$
|| A|| _{(k)}~=\sum _{j=1}^{k}\sigma _{j}(A),\quad (k=1,2, \cdots ,n)$$
and the Schatten $p$-norms defined as
$$
||A|| _{p}=\Big(\sum _{j=1}^{n}(\sigma _{j}(A))^{p}\Big)^{1/p},
$$
$1\leq p< \infty $,
where $\sigma_1\ge \sigma_2\ge\cdots\ge\sigma_n$
are the singular values of $A\in \M_n,$ that is, the eigenvalues of
$|A|$.
The spectral norm (or operator norm) is given by $||A||_\infty =s _{1}(A) =\lim _{p\rightarrow \infty }||A||_p$.

The famous Ky Fan dominance theorem states that a matrix $B$ dominates another matrix $A$ in all UI norms
if and only if it does so in all Ky Fan norms.
The latter set of relations can be written as a weak majorization relation between the vectors of
singular values of $A$ and $B$:
$$
\sigma^\downarrow(A)\prec_w \sigma^\downarrow(B):
\qquad \sum_{j=1}^k \sigma_j(A)\le \sum_{j=1}^k \sigma_j(B), 1\le  k\le n.
$$
For PSD matrices, the above domination relation translates to a weak majorization between
the vectors of eigenvalues: $\lambda^\downarrow(A)\prec_w\lambda^\downarrow(B)$.
Here, $\lambda^\downarrow(A)$ denotes the (real) vector of eigenvalues of $A$ sorted in non-increasing order.

Weyl's monotonicity theorem (\cite{bhatia97}, Corollary III.2.3) states that
$$
\lambda^\downarrow_k(A) \le \lambda^\downarrow_k(A+B), \quad 1\le  k\le n,
$$
for Hermitian $A$ and PSD $B$.

Finally, we refer the reader to Chapter 2 of \cite{kato} for an exposition of a number of
important functional analytic
properties of eigenvalues and corresponding eigenspaces of a Hermitian matrix, which we will need
in the proof of Theorem \ref{theo:delta}.

%%%%%%%%%%%%%%%%%%%%%%%%%%%%%%%%%%%%%%%%%%%%%%%%%%%%%%%%%%%%%%%%%%%%%%%%%%%%%%%%%%%%%%%%%%%%%%%
%%%%%%%%%%%%%%%%%%%%%%%%%%%%%%%%%%%%%%%%%%%%%%%%%%%%%%%%%%%%%%%%%%%%%%%%%%%%%%%%%%%%%%%%%%%%%%%
\section{Comparison of norms $|||f(A)+f(B)|||$ and $|||f(A+B)|||$ \label{sec:main}}
For PD matrices $A,B,$ McCarthy \cite{mccarthy} proved that
$$
||A^r+B^r||_1\le ||(A+B)^r||_1,~~1\le r<\infty
$$
and
$$
||A^r+B^r||_1\ge ||(A+B)^r||_1,~~0\le r\le 1.
$$
Bhatia and Kittaneh \cite{bhatiakittaneh98} proved the above-mentioned inequalities for the operator norm.
There they also proved that
$$
|||A^m+B^m|||\le |||(A+B)^m|||,~~m=1,2,\ldots
$$
for $A,B\ge 0$ and conjectured that if $f$ is operator monotone function on $[0,\infty )$ with $f(0)=0$ then
\begin{equation}\label{eq:1}
|||f(A+B)|||\le |||f(A)+f(B)|||.
\end{equation}
Hiai also posed this conjecture in \cite{hiai}. Ando and Zhan affirmatively settled this conjecture in
\cite{andozhan}. As a corollary they obtained that if $f$ is an increasing function on $[0,\infty )$
with $f(0)=0,~f(\infty )=\infty$ and if the inverse function of $f$ is operator monotone then
\begin{equation}\label{eq:2}
|||f(A+B)|||\ge |||f(A)+f(B)|||.
\end{equation}
Since the inverse function of a non-negative operator convex function on $[0,\infty )$ with $f(0)=0$
is operator monotone \cite{ando88}, we conclude that inequality (\ref{eq:2}) holds for any operator convex
function on $[0,\infty )$ with $f(0)=0$.
In \cite{aujlasilva} it was shown that if the non-negative functions $f,g$ on $[0,\infty )$ satisfy
inequality (\ref{eq:2}) then the functions $f+g,~f\circ g$ and $fg$ also satisfy (\ref{eq:2}).
It was further shown that
any polynomial $p$ with non-negative coefficients and $p(0)=0$ satisfy (\ref{eq:2}).

This prompted the authors
to conjecture in \cite{aujlasilva} that any non-negative convex function on $[0,\infty )$
with $f(0)=0$ should also satisfy (\ref{eq:2}).
Note that such functions must automatically be increasing functions.
Using the fact that a non-negative convex function on $[0,\infty )$
with $f(0)=0$ can be approximated uniformly on a finite interval by a positive linear combination
of angle functions, Kosem settled this conjecture affirmatively
in \cite{kosem}. Later on Bourin and Uchiyama proved (\cite{bourin}; see also \cite{aujlabourin}) that any
non-negative concave function on $[0,\infty)$ (again such functions must be increasing) satisfies (\ref{eq:1}).

It is shown in \cite{aujla,aujlasilva} that
if a non-negative function $f$ satisfies (\ref{eq:1}) then it is concave and
if it satisfies (\ref{eq:2}) then it
is convex with $f(0)=0$. Hence within the set of non-negative $f$ these results
give a full characterisation of all possible $f$ satisfying these inequalities.
This completes our discussion in this section.
%%%%%%%%%%%%%%%%%%%%%%%%%%%%%%%%%%%%%%%%%%%%%%%%%%%%%%%%%%%%%%%%%%%%%%%%%%%%%%%%%%%%%%%%%%%%%%%
%%%%%%%%%%%%%%%%%%%%%%%%%%%%%%%%%%%%%%%%%%%%%%%%%%%%%%%%%%%%%%%%%%%%%%%%%%%%%%%%%%%%%%%%%%%%%%%
\section{Comparison of norms $|||f(A)-f(B)|||$ and $|||f(|A-B|)|||$ \label{sec:mainmin}}
We begin this section with the inequality of Powers and St{\o}rmer \cite{powersstormer},
derived in the course of their work on free states of the canonical anti-commutation relations.
They proved that if $A,B$ are PSD then
$$
||A^{1/2}-B^{1/2}||_{2}^2\le ||A-B||_1.
$$
Kittaneh \cite{kittaneh} generalized this to show that
$$
||A^{1/2}-B^{1/2}||_{2p}^2\le ||A-B||_p
$$
for $1\le p\le \infty$. Note that for any matrix $T$, we have $||T||_{2p}^2=||T^*T||_p,~1\le p\le \infty$,
so this result of Kittaneh can be restated as
$$
||(A-B)^{2}||_{p}\le ||A^2-B^2||_p.
$$
Bhatia \cite{bhatia87} proved  this inequality for all unitarily invariant norms. There he also proved  that
$$
||(A-B)^{2^k}||_{p}\le ||A^{2^k}-B^2{^k}||_p,~~k=1,2,\ldots.
$$
The above inequality when specialized to the $p-$norms gives
$$
||A^{1/m}-B^{1/m}||_{mp}^m\le ||A-B||_p
$$
for all integers $m$ of the form $2^k,~k=1,2,\ldots$, which is an interesting generalisation of the
Powers-St{\o}rmer inequality.

In \cite{birman koplisolo} Birman, Koplienko and Solomyak proved that
$$
||A^r-B^r||_\infty \le ||~|A-B|^r||_\infty ,~~0\le r\le 1
$$
for all  $A,B\ge 0$. Note that the function $f(x)=x^r$ is operator monotone on $[0,\infty )$.
This motivated Kattaneh and Kosaki \cite{kittanehkosaki} to prove that if $f$ is non-negative
operator monotone on $[0,\infty )$  then
$$
||f(A)-f(B)||_\infty \le  f(||A-B||_\infty ) =||f(|A-B|)||_\infty
$$
for all $A,B\ge 0$.
Then Ando \cite{ando88} proved that if $f$ is non-negative operator monotone
on $[0,\infty )$ then
\begin{equation}\label{eq:3}
|||f(A)-f(B)|||\le |||f(|A-B|)|||
\end{equation}
$A,B\ge 0$,  for all unitarily invariant norms. As a corollary to this result, Ando deduced that
the reverse inequality holds for all functions $f$ on $[0,\infty )$ with $f(0)=0$ and
$f(\infty )=\infty $ if the inverse function of $f$ is operator monotone. Since the inverse function
of a non-negative operator convex function on $[0,\infty )$ with $f(0)=0$ is operator monotone
\cite{ando88}, we conclude that if $f$ is operator convex on $[0,\infty )$ with $f(0)=0$ then we have
\begin{equation}\label{eq:4}
|||f(A)-f(B)|||\ge |||f(|A-B|)|||.
\end{equation}

Afterward, Mathias \cite{mathias} proved that the inequality (\ref{eq:3}) holds for any non-negative
matrix monotone  function of order $n$ on $[0,\infty )$.
One may wonder whether, in a similar vein, inequality (\ref{eq:4}) can be proved for a
non-negative increasing matrix convex function $f$ of order $n$ on $[0,\infty )$ with $f(0)=0$.

We have seen that inequality (\ref{eq:1}) holds for non-negative increasing concave functions on $[0,\infty )$
and inequality (\ref{eq:2}) holds for non-negative increasing convex functions on $[0,\infty )$ with $f(0)=0$.
In the same spirit, we consider the question whether inequalities
(\ref{eq:3}) and (\ref{eq:4}) can also be generalized
to non-negative concave and convex functions respectively. We raise and answer several
questions in this direction.

\begin{question}\label{th2a}
For all $A,B,\ge0$, for all UI norms, and for non-negative increasing convex functions
$g$ on $[0,\infty)$ with $g(0)=0$, does the inequality
$||| g(A)-g(B) ||| \ge ||| g(|A-B|) |||$ hold?
\end{question}
The answer to this question is negative, as shown by the following counterexample.
We consider the convex angle function $g(x) = x+(x-1)_+$ and the operator norm.
For the $2\times 2$ PSD matrices
$$
A = \left(\begin{array}{rr}
0.9 &0 \\
0& 0.6
\end{array}
\right),\quad
B = \left(
\begin{array}{rr}
0.8 &0.5 \\
0.5& 0.4
\end{array}
\right),
$$
the eigenvalues of $g(|A-B|)$ are $0.65249$ and $0.35249$,
while those of $g(A)-g(B)$ are $0.65010$ and $-0.48862$. Thus,
$||g(|A-B|)||_\infty =0.65249$, which is larger than
$||g(A)-g(B)||_\infty=0.65010$.
\qed

Under the additional restriction $A\ge B$, the absolute value in the argument
of $g$ in the RHS vanishes, leading to a simplified statement and
a second question, with better hopes for success. Introducing the matrix $\Delta=A-B$,
\begin{question}\label{th2}
For all $B,\Delta\ge0$, for all UI norms, and for non-negative increasing convex functions
$g$ on $[0,\infty)$ with $g(0)=0$, does the inequality
$||| g(B+\Delta)-g(B) ||| \ge ||| g(\Delta) |||$ hold?
\end{question}
This restricted case also turns out to have a negative answer.
Counterexamples, however, were much harder to find, and required a reduction of the problem
based on certain results about a novel majorization-like relation, which we call
dominated majorization.
This will be the subject of Sections \ref{sec:domi} and \ref{sec:apps},
where a number of results of independent interest are proven.

It is also very reasonable to ask:
\begin{question}
For all $B,\Delta\ge0$, for all UI norms, and for non-negative increasing concave functions
$f$ on $[0,\infty)$, does the inequality
$||| f(B+\Delta)-f(B) ||| \le ||| f(\Delta) |||$ hold?
\end{question}

Again, this statement is false, as the following counterexample shows.
Consider the concave angle function $f(x) = \min(x,1) = x-(x-1)_+$,
and the $3\times 3$ PSD matrices
$$
B = \left(
\begin{array}{rrr}
         0.701816  &       0.317887   &      0.198910\\
         0.317887  &       1.014950   &     -0.093826\\
         0.198910  &      -0.093826   &      0.274236
\end{array}
\right)
$$
and
$$
\Delta = \left(
\begin{array}{rrr}
         0.192713 & 0   &      0 \\
         0  &       0.446505 & 0 \\
         0  &       0   &      0.455416
\end{array}
\right).
$$
One gets
$$
||f(\Delta)||_\infty = 0.455416
$$
while
$$
||f(B+\Delta)-f(B)||_\infty = 0.455776.
$$
\qed

Next we consider an even more restricted special case, in which
the inequalities (\ref{eq:3}) and (\ref{eq:4}) finally do hold.
We actually prove that a stronger relationship holds in this special case.
We shall use the notation $\lambda^\downarrow(X)\le \lambda^\downarrow(Y)$ whenever
$\lambda^\downarrow_k(X)\le \lambda^\downarrow_k(Y)$ holds for all $k$.
\begin{theorem}\label{theo:ggc}
For a non-negative,  increasing concave function $g$ on $[0,\infty)$, and
matrices $A,B\ge 0$ such that
$A\ge ||B||_\infty $,
we have
\be\label{eq:ggc}
\lambda^\downarrow(g(A-B)) \ge \lambda^\downarrow(g(A)-g(B)).
\ee
\end{theorem}

An easy corollary is the corresponding statement for non-negative convex functions.
\begin{corollary}\label{prop:gg}
Let $f$ be a non-negative strictly increasing convex function on $[0,\infty )$ with $f(0)=0$.
Let $A,B\ge 0$ be such that $A\ge ||B||_\infty $.
Then
\be\label{eq:gg}
\lambda^\downarrow(f(A-B)) \le \lambda^\downarrow(f(A)-f(B)).
\ee
\end{corollary}
\textit{Proof.}
Let $f=g^{-1}$, with $g$ satisfying the conditions of Theorem \ref{theo:ggc}.
Upon replacing  $A$ by $f(A)$ and $B$ by $f(B)$, the condition $A\ge ||B||_\infty$ is unharmed
as $f$ is monotonous.
Furthermore, (\ref{eq:ggc}) becomes
$$
\lambda^\downarrow(g(f(A)-f(B))) \ge \lambda^\downarrow(A-B).
$$
Applying the function $f$ on both sides does not change the ordering, again because of monotonicity of $f$,
and yields validity of inequality (\ref{eq:gg}).
\qed

These two results obviously imply the corresponding majorization relations, and by Ky Fan dominance,
relations in any UI norm.

%%%%%%%%%%%%%%%%%%%%%%%%%%%%%%%%%%%%%%%%%%%%%%%%%%%%%%%%%%%%%%%%%%%%%%%%%%%%%
\noindent\textit{Proof of Theorem \ref{theo:ggc}.}
W.l.o.g.\ we will assume $||B||_\infty =1$, since any other value can be absorbed in the definition of $g$.

It is immediately clear that if (\ref{eq:ggc}) holds for $g$ that in addition satisfy $g(0)=0$,
then it must also hold without that constraint, i.e.\ for functions $g(x)+c$, with $c\ge 0$.
This is because the additional constant $c$ cancels out in the LHS, while
$\lambda^\downarrow(g(A-B)+c)\ge \lambda^\downarrow(g(A-B))$.

Furthermore, (\ref{eq:ggc}) remains valid when replacing $g(x)$ with $ag(x)$, for $a>0$.
Thus, w.l.o.g.\ we can assume $g(0)=0$ and $g(1)=1$.
Together with concavity of $g$, this implies that, for $0\le x\le 1$, $g(x)\ge x$,
while for $x\ge 1$, the one-sided derivative $g'(x)\le 1$
(since concave functions need not be differentiable everywhere, we have to use the one-sided derivative
$g'(x)=\lim_{t\to 0^+} (g(x+t)-g(x))/t$).

Since $0\le B\le \id$, and for $0\le x\le 1$, $g(x)\ge x$ holds, we have
$g(B)\ge B$, or $-g(B)\le -B$. By Weyl monotonicity, this implies
$\lambda^\downarrow(g(A)-g(B))\le \lambda^\downarrow(g(A)-B)$.
Thus, statement (\ref{eq:ggc}) would be implied by the stronger statement
\be
\lambda^\downarrow(g(A)-B) \le \lambda^\downarrow(g(A-B)). \label{eq:ggc3}
\ee
Now note that the argument of $g$ in the LHS satisfies $A\ge\id$. Thus, in principle, we could
replace $g(x)$ in the LHS by another function $h(x)$ defined as
\be
h(x) = \left\{
\begin{array}{l}
g(x),\mbox{ if }x\ge 1\\
x,\mbox{ otherwise.}
\end{array}
\right.
\ee
If we also do that in the RHS, we get a stronger statement than (\ref{eq:ggc3}).
Indeed, $h(x)\le g(x)$ for $x\ge0$
and $A-B\ge0$, and therefore $h(A-B)\le g(A-B)$ holds. By Weyl monotonicity again, we see that (\ref{eq:ggc3})
is implied by
\be
\lambda^\downarrow(h(A)-B) \le \lambda^\downarrow(h(A-B)). \label{eq:ggc4}
\ee
The importance of this move is that $h(x)$ is still an  increasing and concave function
(because $g'(x)\le 1$ for $x\ge1$),
but now has $h'(x)\le 1$ for $x\ge0$.

Defining $C=A-B$, which is positive semi-definite,
we now have to show the inequality
$$
\lambda^\downarrow_k(h(C+B)-B) \le \lambda^\downarrow_k(h(C)) = h(\lambda^\downarrow_k(C)),
$$
for every $k$.
Fixing $k$, and introducing the shorthand $x_0=\lambda^\downarrow_k(C)$,
we can exploit concavity of $h$ to bound it from above as $h(x)\le a(x-x_0)+h(x_0)$, where
$a=h'(x_0)\le 1$.
Again by Weyl monotonicity, we find
\beas
\lambda^\downarrow_k(h(C+B)-B) &\le& \lambda^\downarrow_k(a(C+B-x_0)+h(x_0)-B) \\
&=& \lambda^\downarrow_k(aC+(a-1)B-ax_0+h(x_0)) \\
&\le& \lambda^\downarrow_k(aC)-ax_0+h(x_0) = h(x_0),
\eeas
where in the second line we could remove the term $(a-1)B$ because it is negative.
This being true for all $k$, we have proved (\ref{eq:ggc4}) and all previous statements that follow from it,
including the statement of the theorem.
\qed
%%%%%%%%%%%%%%%%%%%%%%%%%%%%%%%%%%%%%%%%%%%%%%%%%%%%%%%%%%%%%%%%%%%%%%%%%%%%%%%%%%%%%%%%%%%%%%%
%%%%%%%%%%%%%%%%%%%%%%%%%%%%%%%%%%%%%%%%%%%%%%%%%%%%%%%%%%%%%%%%%%%%%%%%%%%%%%%%%%%%%%%%%%%%%%%
\section{Dominated majorization\label{sec:domi}}
We have already pointed out that inequalities (\ref{eq:1})-(\ref{eq:2}) were proven first for operator
convex or operator concave functions, being extended only afterwards for ordinary convex/concave functions.
Moreover, the proofs for ordinary convex/concave functions actually exploited the corresponding results for
operator convex/concave functions. This may seem somewhat unnatural and it is not unreasonable to ask for
a more direct proof.

In this section we introduce a number of new ideas and techniques which, although they may seem strange and
somewhat contrived
at first, will lead to new, elementary proofs of inequalities (\ref{eq:1})-(\ref{eq:2})
that bypass the Ando-Zhan theorem
and do not require the machinery
of operator monotone and operator convex functions. Secondly, we will use this technique
to try and answer Question 2 raised in the previous section.

%%%%%%%%%%%%%%%%%%%%%%%%%%%%%%%%%%%%%%%%%%%%%%%%%%%%%%%%%%%%%%%%%%%%%%%%%%%%%%%%%%%%%%%%%%%%%%%%%%

Let us consider three Hermitian matrices $A$, $B$ and $C$ and assume that there exists
$a_0>0$
such that the following relation holds for all $a\ge a_0$, and for certain (possibly all) values of $k$:
\be\label{eq:ayb}
\sum_{j=1}^k \lambda_j^\downarrow(aA+B) \le \sum_{j=1}^k \lambda_j^\downarrow(aA+C).
\ee
As it holds for all $a\ge a_0$, it should be possible to simplify this condition.

Subtracting $\sum_{j=1}^k \lambda_j^\downarrow(aA)$ from both sides, and substituting $a=1/t$, we obtain
$$
\frac{1}{t}\sum_{j=1}^k (\lambda_j^\downarrow(A+tB)-\lambda_j^\downarrow(A))
\le
\frac{1}{t}\sum_{j=1}^k (\lambda_j^\downarrow(A+tC)-\lambda_j^\downarrow(A)),
$$
for all $0<t\le t_0 = 1/a_0$.
In the limit of positive $t$ going to 0, this yields a comparison between directional
derivatives of sums of $k$ largest eigenvalues:
\be\label{eq:ddt}
\ddtp\sum_{j=1}^k \lambda_j^\downarrow(A+tB)
\le
\ddtp\sum_{j=1}^k \lambda_j^\downarrow(A+tC).
\ee

Let us introduce the vector $\delta(B;A)$ defined as:
\be\label{eq:defdelta}
\sum_{j=1}^k \delta_j(B;A) := \ddtp\sum_{j=1}^k \lambda_j^\downarrow(A+tB).
\ee
With this notation, relation (\ref{eq:ddt}) becomes
$$
\sum_{j=1}^k \delta_j(B;A) \le \sum_{j=1}^k \delta_j(C;A).
$$
That is, the entries of $\delta(B;A)$ are related via a majorization-like relation
(without the usual rearrangement) to those of $\delta(C;A)$.

To simplify the notations, we will use the symbol $\prec_w$ for this relation:
\be
a \prec_{w} b \Longleftrightarrow \sum_{j=1}^k a_j \le \sum_{j=1}^k b_j,
\ee
and explicitly put rearrangements in the vectors concerned by use of the symbols $\uparrow$ and $\downarrow$.
In that way, we write the classical majorization relation as
$a^\downarrow \prec_w b^\downarrow$.

With these notations
relation (\ref{eq:ddt}) is expressed as
\be\label{eq:ddt2}
\delta(B;A) \prec_{w} \delta(C;A).
\ee
We call this relation \textit{$A$-dominated majorization} or \textit{$A$-majorization} for short.
\begin{definition}
Consider three Hermitian matrices $A$, $B$ and $C$. When the relation (\ref{eq:ddt}) holds,
or equivalently, (\ref{eq:ddt2}), we say that $B$ is $A$-majorized by $C$.
\end{definition}
The argument shown above proves the following:
\begin{proposition}
\label{prop:ainf}
Let $A$, $B$ and $C$ be Hermitian matrices.
If there exists $a_0>0$ such that
$\sum_{j=1}^k \lambda_j^\downarrow(aA+B) \le \sum_{j=1}^k \lambda_j^\downarrow(aA+C)$
holds for all $a\ge a_0$,
then $\delta(B;A) \prec_{w} \delta(C;A)$.
\end{proposition}
%%%%%%%%%%%%%%%%%%%%%%%%%%%%%%%%%%%%%%%%%%%%%%%%%%%%%%%%%%%%%%%%%%%%%%%%%%%
\subsection{Directional derivative of the sum of the $k$-th largest eigenvalues}
It turns out that there is a very simple way to calculate $\delta(B;A)$,
based on an explicit expression of the directional derivative
of the sum of the $k$ largest eigenvalues of a symmetric matrix,
which is well-known in numerical analysis (see \cite{hiriart} and references therein, and \cite{overton}).
The directional derivative of a convex function is defined as follows (\cite{hiriart}, Section 2.2):
\begin{definition}
Let $f(x)$ be a convex function defined on a subset $\mathcal O$ of a Euclidean space $X$.
For any $x\in\mathcal O$, and $d\in X$, the directional derivative of $f$ at $x$ in the direction $d$ is defined as
$$
f'(x,d) = \lim_{t\to 0^+} \frac{f(x+td)-f(x)}{t}.
$$
\end{definition}
It is essential that the limit $t\to 0^+$ is taken because $f$ need not be differentiable.
We will denote this directional derivative by the symbol $\ddtp$.

Consider an $n\times n$ Hermitian matrix $A$, and let its eigenvalues, sorted in non-increasing order, be denoted by
$\lambda_j^\downarrow(A)$, $j=1,2,\ldots,n$. Let its \textit{distinct} eigenvalues, sorted in decreasing order, be denoted
by $\mu_i(A)$, $i=1,2,\ldots,m$ (with $m$ the number of distinct eigenvalues) and the corresponding multiplicities by $r_i$.
Thus $\sum_{i=1}^m r_i=n$.
The sum of the $k$ largest eigenvalues of $A$ can be written in terms of the $\mu_i$ as follows:
writing $k$ as $k=r_1+r_2+\ldots+r_l+s$, where $1\le s\le r_{l+1}$,
$$
\sum_{j=1}^k \lambda^\downarrow_j(A) = \sum_{i=1}^l r_i \mu_i(A) + s \mu_{l+1}(A).
$$
Furthermore, let $P_i$ denote the projector onto the $i$-th eigenspace of $A$, corresponding to eigenvalue $\mu_i(A)$.
Thus, $P_i$ is a matrix of dimensions $r_i\times n$.
The spectral decomposition of $A$ can then be written as
$$
A = \sum_{i=1}^m \mu_i(A) P_i^* P_i.
$$

The following is a reformulation of Corollary 3.9 in \cite{hiriart}, which was proven there for real symmetric matrices.
\begin{proposition}
Let $A$ be a real $n\times n$ symmetric matrix with spectral decomposition
$A = \sum_{i=1}^m \mu_i(A) P_i^* P_i$ and multiplicities $r_i$.
Let $B$ also be a real $n\times n$ symmetric matrix.
With $k$ written as $k=r_1+r_2+\ldots+r_l+s$, where $1\le s\le r_{l+1}$,
the directional derivative of $\sum_{j=1}^k \lambda_j^\downarrow(A)$ in direction $B$ is given by
\be
\ddtp\sum_{j=1}^k \lambda_j^\downarrow(A+tB) =
\sum_{i=1}^l \trace P_i B P_i^* +
\sum_{j=1}^{s} \lambda_j^\downarrow(P_{l+1} B P_{l+1}^*).
\ee
\end{proposition}
Note that, when $s=r_{l+1}$, this formula simplifies to
\be
\ddtp\sum_{j=1}^k \lambda_j^\downarrow(A+tB) = \sum_{i=1}^{l+1} \trace P_i B P_i^*.
\ee

We summarise what we really need to know about this proposition in the following theorem
(quietly extended to the complex case).
\begin{theorem}\label{theo:delta}
Let $A$ and $B$ be Hermitian matrices.
With $\delta(B;A)$ defined by (\ref{eq:defdelta}),
the entries of the vector $\delta(B;A)$ are the diagonal entries of $B$
in a certain basis in which $A$ is diagonal and its diagonal entries
appear sorted in non-increasing order.
When all eigenvalues of $A$ are simple (i.e.\ have multiplicity 1), this basis is just the
eigenbasis of $A$ and does not depend on $B$.
\end{theorem}
An independent proof of this theorem, that also works for complex Hermitian matrices, is presented
in Section \ref{sec:proofdelta}.

The upshot of Theorem \ref{theo:delta} is that there exists a unitary matrix $U$ such that
$U^*AU=\Lambda^\downarrow(A)$ and $\delta(B;A)=\diag(U^*BU)$. In other words, $\delta(B;A)$ is the vector of diagonal
elements of $B$, in a particular basis governed by $A$, and possibly by $B$ too. In the generic case
that all $\lambda_i(A)$ are distinct, $U$ is unique and does not depend on $B$, hence
in that case $\delta(B;A)$ is the vector of diagonal elements of $B$ in the eigenbasis of $A$.
%%%%%%%%%%%%%%%%%%%%%%%%%%%%
\subsection{Dominated majorization for co-diagonal matrices}
Let us now specialise to the case where $A$ and $B$ commute and there is a common basis in which the
diagonal elements of $A$ and $B$ appear in the same, non-increasing order.
We will say that $A$ and $B$ that satisfy this condition are \textit{co-diagonal}.

According to Proposition \ref{prop:ainf}, validity of (\ref{eq:ayb}) for all $a>0$ implies
$A$-majorization, (\ref{eq:ddt2}).
Theorem \ref{theo:delta} now immediately leads to the following proposition, which says that
for co-diagonal $A$ and $B$, validity of (\ref{eq:ayb}) for all $a>0$ is actually equivalent
with $A$-majorization.
\begin{proposition}\label{prop:4b}
For Hermitian $A,B,C$, where $A$ and $B$ are co-diagonal,
the following are equivalent:
\bea
\lambda^\downarrow (aA+B) &\prec_w& \lambda^\downarrow(aA+C), \quad\forall a\ge0 \label{eq:pp1}\\
\delta(B;A) &\prec_{w}& \delta(C;A) \label{eq:pp2} \\
\delta(aA+B;A) &\prec_{w}& \delta(aA+C;A), \quad\forall a\ge0.\label{eq:pp3}
\eea
\end{proposition}
\textit{Proof.}

(\ref{eq:pp1}) implies (\ref{eq:pp2}):\\
If relation (\ref{eq:ayb}) holds for all $a>0$, then it holds for $a$ tending to infinity.
By Proposition \ref{prop:ainf} we then get that $B$ is $A$-majorized by $C$.

(\ref{eq:pp2}) implies (\ref{eq:pp3}): \\
Let us add $a\lambda^\downarrow(A)$ to both sides of (\ref{eq:pp2}).
By Theorem \ref{theo:delta}, $\delta(B;A)$ is the vector of diagonal elements of $B$,
in a basis in which $A$ is diagonal and the eigenvalues of $A$
appear sorted in non-increasing order.
Thus, $\forall a>0$, $\delta(B;A) + a\lambda^\downarrow(A) = \delta(B+aA;A)$.
The same holds for $C$.

(\ref{eq:pp3}) implies (\ref{eq:pp1}):\\
By the co-diagonality of $A$ and $B$, $aA+B$ is diagonal in any basis in which $A$ is diagonal.
Hence, the LHS of (\ref{eq:pp3}) is equal to $\lambda^\downarrow(a A+B)$.
By Schur's majorization theorem, the RHS of (\ref{eq:pp3}) is majorized by $\lambda^\downarrow(aA + C)$.
\qed
%%%%%%%%%%%%%%%%%%%%%%%%%%%%%%%%%%%%%%%%%%%%%%%%%%%%%%%%%%%%%%%%%%%%%%%%%%%%%%%%%%%%%%%%%%%%%
\section{Applications of dominated majorization\label{sec:apps}}
In this section we first use Proposition \ref{prop:4b},
to give a new, elementary proof of inequality (\ref{eq:1}) for non-negative concave functions
(which readily implies validity of inequality (\ref{eq:2}) for non-negative convex functions),
that does not rely on the Ando-Zhan inequality
for operator concave functions, nor on the theory of operator monotone functions.

Then, we answer Question 2 in the negative by exhibiting a counterexample.
Here, too, Proposition \ref{prop:4b} was instrumental.
%%%%%%%%%%%%%%%%%%%%%%%%%%%%%%%%%%%%%%%%%%%%%%%%%%%%%%%%%%%%%%%%%%%%%%%%%%%%%
\subsection{A new proof of inequality (\ref{eq:1}) for non-negative concave functions}
We want to prove that
$$
|||f(A+B)|||\le |||f(A)+f(B)|||
$$
holds for all non-negative concave functions $f(x)$.
Therefore, it should hold in particular for all functions $f(x)=b+ax+f_0(x)$,
where $f_0$ is non-negative concave with $f_0(0)=0$ and $f_0'(+\infty)=0$,
and for all $a,b\ge0$.
Inserting this in the eigenvalue-majorization form of inequality (\ref{eq:1}), we get the
majorization relation
$$
\lambda^\downarrow(b+a(A+B)+f_0(A+B))\prec_w \lambda^\downarrow(2b+ a(A+B)+f_0(A)+f_0(B)),
$$
for $A,B\ge0$.
Clearly, this is strongest for $b=0$.
Proposition \ref{prop:4b} then immediately yields the equivalent form
$$
\delta(f(A+B);A+B)\prec_{w} \delta(f(A)+f(B);A+B),
$$
for all non-negative concave functions $f$ (recall that such functions are non-decreasing) with $f(0)=0$.

An interesting aspect of this form is that, unlike $\lambda$, $\delta$ is linear in its first argument.
Our proof of the equivalent form, stated as Proposition \ref{prop:bourins}
below, crucially depends on this property.

\begin{proposition}\label{prop:bourins}
For positive semidefinite $A$ and $B$, and $f$ a non-negative concave function with $f(0)=0$,
\be\label{eq:bourins}
\delta(f(A+B);A+B)\prec_{w} \delta(f(A)+f(B);A+B).
\ee
\end{proposition}
\textit{Proof.}
Any non-negative concave function $f$ can be uniformly approximated as a positive linear combination of
angle functions $x\mapsto x-(x-t)_+$.
By linearity of $\delta$, inequality (\ref{eq:bourins}) follows if it holds for any such angle function, i.e.
$$
\delta(A+B-(A+B-t)_+;A+B)\prec_{w} \delta(A-(A-t)_+ + B-(B-t)_+;A+B),
$$
which, again by linearity, simplifies to
$$
\delta((A-t)_+ +(B-t)_+;A+B)\prec_w \delta((A+B-t)_+;A+B).
$$
In fact, for angle functions the latter inequality even holds with rearrangement, and we shall prove
$$
\delta^\downarrow((A-t)_+ +(B-t)_+;A+B)\prec_w \delta^\downarrow((A+B-t)_+;A+B),
$$
for all $t\ge0$.
Letting $\tr(x)$ denote the sum $\sum_{i=1}^n x_i$ of $x=(x_1,\ldots,x_n)$,
this relation can be expressed in a well-known way as
$$
\tr(\delta((A-t)_+ +(B-t)_+;A+B)-s)_+ \le \tr(\delta((A+B-t)_+;A+B)-s)_+,
$$
for all $s$ (and $t\ge0$). Since both vectors $\delta$ are non-negative it suffices to consider the case $s\ge0$.
In the eigenbasis of $A+B$, $A+B$ itself is of course diagonal, hence
the RHS simplifies to $\trace(A+B-(s+t))_+$.

Now we introduce the variable $u=s+t$.
The last inequality has to be valid for all values of $s$ and $t$, thus if we keep the value of $u$ fixed,
the inequality has to remain true if we maximise the LHS over all values of $t$
in the range $[0,u]$ (and set $s=u-t$).
That is,
\bea\label{eq:tr1}
\lefteqn{
\max_{0\le t\le u} \tr(\delta((A-t)_+ +(B-t)_+;A+B)-u+t)_+} \nonumber \\
&\le& \trace(A+B-u)_+.
\eea

The next important consequence of the simple behaviour of $\delta$
is that the function $t\mapsto F(t):=\tr(\delta((A-t)_+ +(B-t)_+;A+B)-u+t)_+$
is convex.
Note first that the positive part function is convex and increasing. Applying this to its
outer appearance in the definition of $F$, the required convexity of $F(t)$ follows if, for any $i$,
$\delta((A-t)_+ +(B-t)_+;A+B)_i-u+t$ is itself a convex function of $t$. This function
can be written as $((A-t)_+)_{ii} +((B-t)_+)_{ii}-u+t$, in the eigenbasis of $A+B$.
Hence, convexity follows from the convexity of $t\mapsto (\psi,(A-t)_+\psi)$, for any vector $\psi$,
and to see the latter, just consider this quantity in the eigenbasis of $A$ and see that it can be written as
$\sum_{j=1}^n (\lambda_j(A)-t)_+|\psi_j|^2$, which is a positive linear combination of angle functions
and, therefore, convex.

The convexity of $F(t)$ now implies the simple fact that the maximum in the LHS of (\ref{eq:tr1})
$\max_{0\le t\le u} \tr(\delta((A-t)_+ +(B-t)_+;A+B)-u+t)_+$
is achieved in one of the extreme points, either in $t=0$ or in $t=u$. Noting that $A$ and $B$
are positive semidefinite,
the value achieved in $t=0$ is $\tr(\delta(A +B;A+B)-u)_+$, which is identical to
the RHS in (\ref{eq:tr1}). It therefore only remains to show that the value in $t=u$ is also
bounded above by the RHS. Using the fact
that the function $\tr\delta(X;Y)$ is always equal to $\trace X$,
this amounts to the inequality
\be\label{eq:tru}
\trace(A-u)_+ +\trace(B-u)_+ \le \trace(A+B-u)_+.
\ee
Here, the outer appearance of the positive part function in the LHS has been removed
because its argument is always positive semidefinite.

To prove inequality (\ref{eq:tru}), recall the norm inequality
$$
||| A\oplus B||| \le |||(|A|+|B|)\oplus 0|||,
$$
valid for any unitarily invariant norm (\cite{bhatia97}, Theorem IV.2.13). In particular, it holds for
the Ky Fan norms, and for PSD $A$ and $B$ can be written as the eigenvalue majorization
$$
\lambda^\downarrow\left(A\oplus B\right) \prec_w \lambda^\downarrow\left((A+B)\oplus 0\right).
$$
Thus, for all $u\ge0$ (again, by non-negativity of $A$ and $B$, it suffices to consider $u\ge0$),
$$
\trace\left(\twomat{A}{0}{0}{B}-u\right)_+ \le
\trace\left(\twomat{A+B}{\phantom{A}0}{0}{\phantom{A}0}-u\right)_+,
$$
which is nothing but inequality (\ref{eq:tru}), reformulated in terms of $2\times 2$ block matrices.
This ends the proof of the proposition.
\qed

One might still object that our proof is not really elementary, relying as it is on Proposition \ref{prop:4b}
and the theory behind it.
Strictly speaking, though, Proposition \ref{prop:4b} is not needed in the proof, and only provided the intuition
to try and prove the equivalent form (\ref{eq:bourins}).
Indeed, validity of inequality (\ref{eq:1}) follows immediately from Proposition \ref{prop:bourins} by
combining it with Schur's majorization theorem:
\beas
\lambda^\downarrow(f(A+B)) &=& \delta(f(A+B);A+B) \\
&\prec_{w}& \delta(f(A)+f(B);A+B) \\
&\prec_w& \lambda^\downarrow(f(A)+f(B)).
\eeas

As already shown by Ando and Zhan \cite{andozhan}, validity of inequality (\ref{eq:1}) for a given
non-negative increasing concave
function $f$ implies inequality (\ref{eq:2}) for the inverse function $g=f^{-1}$.
Hence, in combination with our proof of inequality (\ref{eq:1}), this also yields an elementary
proof of inequality (\ref{eq:2}) for non-negative convex functions $g(x)$ with $g(0)=0$,
This was first proven independently from (\ref{eq:1}) by Kosem,
by appealing to the corresponding inequality for operator convex functions.

For completeness, we repeat the short Ando-Zhan argument here.

\textit{Proof of inequality (\ref{eq:2}) for non-negative convex functions.}
Let $g(x)$ be a non-negative convex function with $g(0)=0$. Thus $g$ is increasing.
In particular, $g$ applied to vectors is strongly isotone \cite{bhatia97}.

Let $f(x)$ be its inverse function, $f=g^{-1}$;
thus $f(x)$ is a non-negative increasing concave function with $f(0)=0$.
For such $f$, we have (inequality (\ref{eq:1}))
$$
\lambda^\downarrow(f(A+B)) \prec_w \lambda^\downarrow(f(A)+f(B)).
$$
Since $g(x)$ is strongly isotone, applying $g$ on both sides preserves weak majorization:
$$
g(\lambda^\downarrow(f(A+B))) \prec_w g(\lambda^\downarrow(f(A)+f(B))).
$$
This simplifies, by monotonicity of $g$, to
$$
\lambda^\downarrow(A+B) \prec_w \lambda^\downarrow(g(f(A)+f(B))).
$$
Substituting $g(A)$ for $A$ and $g(B)$ for $B$ then yields inequality (\ref{eq:2}).
\qed
%%%%%%%%%%%%%%%%%%%%%%%%%%%%%%%%%%%%%%%%%%%%%%%%%%%%%%%%%%%%%%%%%%%%%%%%%%%%%%%%%%%%%%%
\subsection{Counterexample to Question \ref{th2}\label{sec:counter}}
To answer Question \ref{th2}, we  will first disregard the absolute values and consider the property
that a convex function $f$ satisfies
\be\label{eq:star1}
\lambda(f(\Delta)) \prec_w \lambda(f(B+\Delta)-f(B))
\ee
for all PSD $B$ and $\Delta$,
which is equivalent to the statement
\be\label{eq:star2}
\lambda(f(A-B)) \prec_w \lambda(f(A)-f(B))
\ee
for all $A\ge B\ge0$.

Although it is by no means obvious at this point, when $\Delta>0$ strictly, Question \ref{th2} is equivalent
to validity of (\ref{eq:star1}) for all stated functions. While it is obvious that (\ref{eq:star1}) implies
$||| g(B+\Delta)-g(B) ||| \ge ||| g(\Delta) |||$, the opposite is not necessarily true because of the
absolute value implicit in the definition of the norm. Nevertheless, it will turn out that a counterexample to
(\ref{eq:star1}) for some function $g$ will indirectly yield a counterexample to Question \ref{th2}
for \textit{some other function}
$\tilde{g}(x)=g(x)+\alpha x$, with $\alpha>0$ large enough, provided $\Delta>0$ holds strictly.
Here, $\alpha$ must be large enough to make $\tilde{g}(B+\Delta)-\tilde{g}(B)=g(B+\Delta)-g(B)+\alpha\Delta$
positive semidefinite,
in which case the absolute value signs can be left out. This will all be made clear below.

The monotone convex angle functions $x\mapsto ax+(x-1)_+$ ($a\ge0$) already have proven
their valour as a testing ground for
similar statements, in Section \ref{sec:main}.
Numerical experiments using angle functions for inequality (\ref{eq:star1}) did not directly lead
to any counterexamples, however. This temporarily increased our belief that the inequality might actually hold,
and led us to investigate, as an initial step towards a `proof', whether the inequality
$$
\sum_{j=1}^k \lambda_j^\downarrow(aY+B) \le \sum_{j=1}^k \lambda_j^\downarrow(aY+C)
$$
might be true for all $a\ge0$, where $B=f(Y)$ and $C=f(X+Y)-f(X)$, and $f(x)=(x-1)_+$.

\bigskip

If the answer to Question \ref{th2} is to be affirmative, it should at least
hold for all angle functions $f(x)=ax+b(x-x_0)_+$.
By Proposition \ref{prop:4b} this is equivalent to the statement
$$
\delta((Y-\id)_+;Y) \prec_{w} \delta((X+Y-\id)_+ - (X-\id)_+;Y).
$$
Consider the $3\times 3$ PSD matrices
$$
X = \left(
\begin{array}{rrr}
         0.35614  &       -0.053243   &      0.10116\\
         -0.053243  &       0.87456   &     0.40559\\
         0.10116  &      0.40559   &      0.82474
\end{array}
\right)
$$
and
$$
Y = \left(
\begin{array}{rrr}
         0.53642  &                       0   &                      0\\
                         0  &       0.42018   &                      0\\
                         0  &                       0   &      0.094866
\end{array}
\right).
$$
The eigenbasis of $Y$ is therefore the standard basis.
Then $\delta((Y-\id)_+;Y)=(0,0,0)$
and
$$
(X+Y-\id)_+ - (X-\id)_+ = \left(
\begin{array}{rrr}
-0.00018194  & 0.00052449 &  -0.0016345\\
   0.00052449 &      0.2573 &     0.12368\\
   -0.0016345  &    0.12368 &        0.04
\end{array}
\right)
$$
so that $\delta((X+Y-\id)_+ - (X-\id)_+;Y)=(-0.00018194, 0.2573, 0.04)$.
The first entry is negative, violating the majorization relation.

Now, as mentioned above, this counterexample immediately yields a counterexample to Question \ref{th2}.
Consider thereto the function $f(x)=\alpha x+(x-1)_+$ with $\alpha=1$, say.
Then the LHS of the inequality becomes
$\delta(Y+(Y-\id)_+;Y) = (0.53642, 0.42018, 0.094866)$ and the RHS
$\delta(Y+(X+Y-\id)_+ - (X-\id)_+;Y) = (0.53624,0.67748,0.13487)$,
again violating the inequality. Since
$Y+(X+Y-\id)_+ - (X-\id)_+$ is a positive definite matrix (as can be checked numerically),
it is unchanged by putting in the required absolute value signs.

Even more explicitly, consider the function
$g(x)= 101 x + (x-1)_+$. Then
$\lambda^\downarrow(g(X+Y)-g(X)) = (54.17824, 42.69595, 9.621004)$ while
$\lambda^\downarrow(g(Y))        = (54.17842, 42.43818, 9.581466)$.
This clearly violates the eigenvalue majorization relation of Question \ref{th2}, \textit{with} absolute
value signs, because of the positivity of $g(X+Y)-g(X)$.

%%%%%%%%%%%%%%%%%%%%%%%%%%%%%%%%%%%%%%%%%%%%%%%%%%%%%%%%%%%%%%%%%%%%%%%%%%%%%%%%%%%%%%%%%%
%\pagebreak
\section{Proof of Theorem \ref{theo:delta}\label{sec:proofdelta}}
In this section, we give a self-contained proof of Theorem \ref{theo:delta} that does not rely on
the methods of convex analysis and is also valid for complex Hermitian matrices, not only real-symmetric ones.
For convenience, we reformulate the statement of the theorem here.

Define a \textit{proper eigenbasis} of a Hermitian matrix $A$ as an orthonormal basis in which 
$A$ is diagonal and its diagonal entries are the eigenvalues of $A$
sorted in non-increasing order.

\begin{quote}\textbf{Theorem 2'.}
Let $A$ and $B$ be Hermitian matrices.
With $\delta(B;A)$ defined via equation (\ref{eq:defdelta}),
the entries of the vector $\delta(B;A)$ are the diagonal entries of $B$
in some proper eigenbasis of $A$.
When all eigenvalues of $A$ are simple (i.e.\ have multiplicity 1), this proper eigenbasis is unique;
otherwise the required one depends on $B$.
\end{quote}

We need a number of definitions first, and recall some basic facts about the perturbation theory of 
eigenvalue decompositions
(see, e.g.\ \cite{kato}, Chapter 2, Section 1).

Consider the matrix-valued function $z\mapsto A+zB$, $z\in\C$,
with $A$ and $B$ the $n\times n$ Hermitian matrices of the theorem.
It is well-known that the roots of the characteristic function of $A+zB$ are
analytic functions of $z$ with only algebraic singularities.
This means that the number $m$ of (distinct) eigenvalues of $A+zB$ is a constant of $z$,
with the exception of a number of special values of $z$, which will be
called exceptional points. If $m<n$, we say that $A+zB$ is permanently degenerate.
In the exceptional points some of the eigenvalues may coincide; this is called an
accidental degeneracy.

In the following we will consider a simply-connected subdomain $D$ of the complex plane $\C$
containing no exceptional points, and such that the intersection of $D$ with the real axis
is the interval $(0,t_0)$, with $t_0>0$. The closure of $D$ is denoted $\overline{D}$, and its intersection
with the real axis is $[0,t_0]$.

We can write the (possibly multiple) eigenvalues of $A+zB$, $z\in \overline{D}$, as holomorphic functions
$\lambda_1(z), \lambda_2(z),\ldots, \lambda_n(z)$. For $z=t\in\R$, these eigenvalues are real
and can be sorted. Sorted in non-increasing order they will be denoted as
$\lambda_1^\downarrow(t) \ge \lambda_2^\downarrow(t) \ge \ldots\ge \lambda_n^\downarrow(t)$, $0\le t\le t_0$.

Furthermore, we can write the \textit{distinct} eigenvalues of $A+zB$, $z\in D$, as a fixed number
of holomorphic functions
$\mu_1(z), \mu_2(z),\ldots, \mu_m(z)$. We will number them such that
$\mu_1(t)>\mu_2(t)>\ldots>\mu_m(t)$ holds for $t\in(0,t_0)$ (or $t\in[0,t_0)$).
We denote the multiplicity of $\mu_i$ by $r_i$. Thus, $n=\sum_{i=1}^{m} r_i$.

The projector on the eigenspace of $A+zB$ corresponding to $\mu_i(z)$ will be denoted
by the function $\cP_i(z)$, $z\in D$, and is called the eigenprojection for $\mu_i(z)$.
This function is holomorphic on $D$ \cite{kato}.

If $z=0$ is not an exceptional point, the distinct eigenvalues of $A$ are equal to the limiting values
$\mu_i(0)$, and the corresponding eigenprojections coincide with the $\cP_i(0)$.

If $z=0$ is an exceptional point then an accidental degeneracy occurs and 
$A$ has less than $m$ distinct eigenvalues.
Each of these eigenvalues may split into several $\mu_i(t)$;
that is, $\lim_{t\to0}\mu_i(t)=\lambda$ for several (contiguous) values of $i$, say $i=i_1,\ldots,i_2$, 
where $\lambda$ is a certain eigenvalue of $A$.
In that case, the eigenprojection for $\lambda$ of $A$ coincides with the sum 
$\lim_{t\to0}\sum_{j=i_1}^{i_2} \cP_{j}(t)$; i.e.\ the
$\lim_{t\to0} \cP_{j}(t)$ separately are not themselves eigenprojectors of $A$.

\bigskip

Let $k$ be an integer such that there exists an $l$ for which $k=r_1+r_2+\ldots+r_l$; we shall say that
such a $k$ is an \textit{entire sum} of the multiplicities $r_i$.
For such values of $k$, we define the projector $\cP_{(k)}(z)$ as the sum of eigenprojectors
$$
\cP_{(k)}(z) = \cP_1(z)+\cP_2(z)+\ldots+\cP_l(z).
$$
For $z=t$ real, this is the projector on the subspace spanned by the 
eigenvectors of the $k$ largest eigenvalues (counting
multiplicities) of $A+tB$.
Since the $\cP_i(z)$ are holomorphic functions on $D$, so is $\cP_{(k)}(z)$.
By continuity of the eigenvalues $\lambda_k(z)$, we have for any such $k$,
$$
\sum_{j=0}^k \lambda^\downarrow_k(A) = \lim_{t\to 0^+}\sum_{j=0}^k \lambda^\downarrow_k(t)
=\trace[\lim_{t\to 0^+}\cP_{(k)}(t)\,\, A].
$$

If $k$ cannot be written in this way, i.e.\ $k=r_1+r_2+\ldots+r_l+s$ with $s$ a `remainder'
satisfying $0<s<r_{l+1}$,
we cannot uniquely define $\cP_{(k)}(z)$, because there is an infinity of $s$-dimensional subspaces
in the eigenspace for $\mu_{l+1}$.
Hence, we will only define $\cP_{(k)}(z)$ for $k$ that are entire sums of $r_i$.

Finally, we define the projectors $\cP_{(k)}$.
If $z=0$ is not an exceptional point, and $k$ is an entire sum of multiplicities $r_i$, then
$\cP_{(k)}(0)$ is defined, and we define $\cP_{(k)}:=\cP_{(k)}(0)$.
If $z=0$ is an exceptional point then $A+zB$ has an accidental degeneracy at $z=0$. Even if $k$ is an entire
sum of multiplicities $r_i$ of $A+zB$, it need not be an entire sum of multiplicities of $A$.
Hence, in that case $\cP_{(k)}(t)$ is only defined for $t\in(0,t_0)$ (with $t_0>0$) but
not for $t=0$.
We will then define $\cP_{(k)}$ as the limiting value
$$
\cP_{(k)} := \lim_{t\to 0^+} \cP_{(k)}(t).
$$
For all other values of $k$, $\cP_{(k)}$ will not be defined.

\bigskip

\begin{lemma}
If $k$ is such that $\cP_{(k)}$ is defined (directly in $t=0$ or via the limit $t\to 0^+$), then
$$
\sum_{j=1}^k \delta_j(B;A) = \trace B\cP_{(k)}.
$$
\end{lemma}
\textit{Proof.}
Consider the variational characterization of the sum of the $k$ largest eigenvalues of
a Hermitian matrix $H$:
$$
\sum_{j=1}^k \lambda_j^\downarrow(H) = \max_{Q} \trace[H\,Q],
$$
where $Q$ runs over all rank-$k$ projectors.
If $k$ is such that $\cP_{(k)}(H)$ exists (taking the potential degeneracies of $H$ into account)
then $Q=\cP_{(k)}(H)$ achieves the maximum,
i.e.\ $\max_{Q} \trace[H\,Q] = \trace[H\,\cP_{(k)}(H)]$.

We have, in particular, that $\cP_{(k)}(t):=\cP_{(k)}(A+tB)$ (if it exists) achieves the maximum for $H=A+tB$.
More precisely, for any $t$ in the open interval $(0,t_0)$, 
the function $u\mapsto \trace [(A+tB)\cP_{(k)}(u)]$ achieves its maximum over $(0,t_0)$
in the \textit{interior} point $u=t$.
Since $\cP_{(k)}(t)$ is holomorphic, this function is differentiable, hence
this maximum must be a stationary point.
Thus
$$
\frac{\partial}{\partial u}\Bigg|_{u\to t}\trace[(A+tB)\,\cP_{(k)}(u)]=0,
$$
i.e.
$$
\trace[(A+tB)\,\ddt\cP_{(k)}(t)]=0.
$$

This implies
\beas
\ddt\sum_{j=1}^k \lambda_j^\downarrow(t)
&=& \ddt\trace[(A+tB)\,\cP_{(k)}(t)] \\
&=& \trace[(A+tB)\, \ddt\cP_{(k)}(t)] + \trace[B\,\cP_{(k)}(t)] \\
&=& \trace[B\,\cP_{(k)}(t)].
\eeas
In particular,
$$
\sum_{j=1}^k \delta_j(A;B)
= \ddtp\sum_{j=1}^k \lambda_j^\downarrow(t)
= \lim_{t\to 0^+}\trace[B\,\cP_{(k)}(t)]
= \trace[B\,\cP_{(k)}].
$$
\qed

\bigskip

We are now in the position to prove Theorem \ref{theo:delta}.
Let's first consider the simplest case when $A$ is not degenerate, i.e.\ all eigenvalues of $A+zB$
are simple for $z\in\overline{D}$.
In that case $\cP_{(k)}(z)$ is always defined for all $k$ and all $z\in\overline{D}$, and, hence,
$\cP_{(k)}$ is defined as $\cP_{(k)}(0)=\sum_{j=1}^k \cP_j(0)$.
There is a unique unitary matrix $U$ such that $UAU^*=\diag(\lambda^\downarrow(A))$,
and in this basis the projector $\cP_j$ is expressed as $e^{jj}$.
Hence, by the lemma we have that $\sum_{j=1}^k \delta_j(B;A) = \trace B\cP_{(k)} = \sum_{j=1}^k B_{jj}$,
where the $B_{jj}$ are the diagonal elements of $B$ expressed in that same basis.
Therefore, for all $j$, $\delta_j(B;A) = B_{jj}$.

If $A$ is degenerate, there is no unique eigenbasis of $A$.
However, the lemma only requires us to deal with the limits $\lim_{t\to 0^+}\cP_{(k)}(t)$.
If the degeneracy of $A$ is lifted completely in $A+zB$, i.e.\ all eigenvalues of $A$ split into simple eigenvalues,
then all $\cP_j(t)$ are rank-1 projectors and $\cP_{(k)}(t)=\sum_{j=1}^k \cP_j(t)$.
Furthermore, letting $i_1$ and $i_2$ be any pair of indices such that an eigenvalue of $A$ splits into 
the eigenvalues $\mu_{i_1},\ldots,\mu_{i_2}$ of
$A+zB$, we have that $\lim_{t\to0}\sum_{j=i_1}^{i_2} \cP_{j}(t)$ is an eigenprojector of $A$.
Therefore, there exists a unique proper eigenbasis of $A$ (determined by $B$)
in which $\lim_{t\to0}\cP_{j}(t)=e^{jj}$, the elementary matrix with a $1$ in position $(j,j)$ and zeroes elsewhere.
Again we find that,
for all $j$, $\delta_j(B;A) = B_{jj}$ in that proper eigenbasis.

The most complicated case arises
when $A+zB$ is permanently degenerate, i.e.\ the degeneracies are not lifted completely, as some
eigenvalues of $A$ may split into still degenerate eigenvalues $\mu_i$ of $A+zB$, with multiplicities $r_i$.
Then the projectors $\cP_i(z)$ have rank $r_i$,
and the $\cP_{(k)}(t)$ are only defined
when $k$ is an entire sum of the multiplicities $r_i$.
There still exists a proper eigenbasis of $A$
in which the projectors $\lim_{t\to0}\cP_{j}(t)$ are diagonal, now of the form $0\oplus \id_{r_j}\oplus 0$, but it
is no longer unique; we will exploit exactly this freedom to deal with $k$ that are not entire sums.

If $k$ is not an entire sum of $r_i$,
we have $k=r_1+r_2+\ldots+r_l+s$, with $s$ the remainder term, satisfying
$1\le s< r_{l+1}$. 
We first write $k$ as an interpolated value between two entire sums as follows:
\beas
k &=& \frac{s}{r_{l+1}}(r_1+r_2+\ldots+r_{l+1}) + (1-\frac{s}{r_{l+1}})(r_1+r_2+\ldots+r_{l}) \\
&=& \alpha k_+ + (1-\alpha) k_-.
\eeas
Here we defined $\alpha=s/r_{l+1}$, and the two entire sums
$k_- = r_1+r_2+\ldots+r_{l}$ and $k_+ = r_1+r_2+\ldots+r_{l+1}$.
We can express $\sum_{j=1}^k \lambda^\downarrow_j$
as a linear interpolation between $\sum_{j=1}^{k_-} \lambda^\downarrow_j$
and $\sum_{j=1}^{k_+} \lambda^\downarrow_j$:
\beas
\sum_{j=1}^k \lambda^\downarrow_j(t)
&=& \sum_{i=1}^{l} r_i \mu_i(t) + s\mu_{l+1}(t) \\
&=& \trace[(A+tB)\,(\cP_1(t)+\ldots+\cP_l(t)+\frac{s}{r_{l+1}}\cP_{l+1}(t))] \\
&=& \trace[(A+tB)\,(\alpha \cP_{(k_+)}(t) + (1-\alpha)\cP_{(k_-)}(t))]\\
&=& \alpha \trace[(A+tB)\,\cP_{(k_+)}(t)] + (1-\alpha)\trace[(A+tB)\,\cP_{(k_-)}(t))].
\eeas
Applying the Lemma to both terms,
we obtain
\be
\sum_{j=1}^k\delta_j(B;A) 
= \trace[B(\alpha \cP_{(k_+)} + (1-\alpha)\cP_{(k_-)})]
= \sum_{i=1}^{l} \trace B\cP_i +\alpha\trace B\cP_{l+1}. \label{eq:interpol}
\ee
Again, to deal with eigenvalue splitting at $z=0$,
each of the $\cP_i$ corresponds to the
limit $\lim_{t\to 0^+}\cP_i(t)$.

Let us consider a partitioning of $B$ in an eigenbasis of $A+zB$ mentioned before, in which the $\cP_i(z)$ 
appear in the form $0\oplus \id_{r_i}\oplus0$. That is,
in $B$ we can single out blocks on its diagonal, each of which corresponds to an eigenspace of $A+zB$; 
Then $\trace B\cP_{i}(z)$ is the sum of all $r_i$ diagonal elements of the $i$-th block of $B$.

The degeneracy of the eigenvalues $\mu_i(z)$ means that this eigenbasis is still not unique
and is determined up to `local' rotations within each of the eigenspaces.
We can use this freedom to make the diagonal elements of $B$ equal within each block. 
This allows us to get rid of $\alpha$ in (\ref{eq:interpol}).
Indeed, as $\alpha = s/r_{l+1}$ and $\trace B\cP_{l+1}$ is the sum of all $r_{l+1}$ diagonal elements of
the $(l+1)$-th block of $B$, then if all these diagonal elements are equal,
$\alpha\trace B\cP_{l+1}(z)$ is equal to the sum of the first $s$ diagonal elements of $B$ in that block.

Wrapping up we find that $\sum_{i=1}^{l} \trace B\cP_i +\alpha\trace B\cP_{l+1}$ equals the sum of the
first $r_1+r_2+\ldots+r_l+s = k$ diagonal elements of $B$ in the chosen eigenbasis.
Taking the limit $z=t\to 0$, we finally obtain
that, again, there is a proper eigenbasis of $A$ in which
$$
 \sum_{j=1}^k\delta_j(B;A)  = \sum_{j=1}^{k} B_{jj},
$$
and hence $\delta_j(B;A)=B_{jj}$.
\qed

%%%%%%%%%%%%%%%%%%%%%%%%%%%%%%%%%%%%%%%%%%%%%%%%%%%%%%%%%%%%%%%%%%%%%%%%%%%%%%%%%%%%%%%%%%%%%%%%%%%%%%%%%%

\bigskip

KA acknowledges the hospitality of the Institut Mittag-Leffler, Djursholm (Sweden), where the final stages of the
work have been done.
We thank an anonymous referee for a variety of detailed comments, which helped to improve the exposition considerably.

%------------------------------------------------------------- BIBLIOGRAPHY

%%%%%%%%%%%%%%%%%%%%%%%%%%%%%%%%%%%%%%%%%%%%%%%%%%%%%%%%%%%%%%%%%%%
\end{document}